\documentstyle[12pt]{article}
\begin{document}
\begin{center}
{\large \bf SOME REMARKS ON PRODUCING HOPF ALGEBRAS}
\end{center}
\begin{center}
{\sc A. A. Vladimirov}
\end{center}
\begin{center}
{\em Bogoliubov Laboratory of Theoretical Physics,
Joint Institute for Nuclear Research, \\
Dubna, Moscow region 141980, Russia}
\end{center}

\vspace{.5cm}

{\small We report some observations concerning two well-known
approaches to construction of quantum groups. Thus, starting from a
bialgebra of inhomogeneous type and imposing quadratic,
cubic or quartic commutation relations on a subset of its
generators we come, in each case, to a $q$-deformed universal enveloping
algebra of a certain simple Lie algebra. An interesting correlation
between the order of initial commutation relations and the Cartan
matrix of the resulting algebra is observed. Another example
demonstrates that the bialgebra structure of $sl_q(2)$ can be
completely determined
by requiring the $q$-oscillator algebra to be its covariant
comodule, in analogy with Manin's approach to define $SL_q(2)$ as a
symmetry algebra of the bosonic and fermionic quantum planes.}

\vspace{1cm}

A lot of recipes are known to produce new examples
of quantum groups (quasitriangular Hopf algebras). The aim of the
present note is to demonstrate some interesting results which
can be obtained within the two of these approaches. Although the
Hopf algebras thus obtained proved to be not new, the way of
producing them looks instructive and fruitful.

Our first example deals with the Hopf algebras generated by a
constant matrix solution $R$ of the Yang-Baxter equations.
The appropriate technique has been elaborated in~\cite{FRT,Ma-rev}.
Here we actually use its modification described in~\cite{Lu,Vl2}.
Let us consider a bialgebra of inhomogeneous type~\cite{SWW} with
generators $\{1,t_{j}^{i},u_{j}^{i},E_i,F^i\}$ which form matrices
$T,U$, a row $E$ and a column $F$, respectively. Its multiplication
relations are
\begin{equation}
\label{1} \begin{array}{ll}
R_{12}\,T_1\,T_2=T_2\,T_1\,R_{12}\,, & E_1\,T_2=T_2\,E_1\,
R_{12}\,, \\
R_{12}\,U_1\,U_2=U_2\,U_1\,R_{12}\,, & F_2\,U_1=
R_{12}\,U_1\,F_2\,, \\
R_{12}\,U_1\,T_2=T_2\,U_1\,R_{12}\,, & T_2\,F_1=
R_{12}\,F_1\,T_2\,, \\
E\,F-F\,E=T-U\,, & U_1\,E_2=E_2\,U_1\,R_{12}\,,
\end{array}
\end{equation}
the coalgebra structure is
\begin{equation}
\label{2}  \begin{array}{ll}
\Delta(T)=T\otimes T\,, & \Delta(E)=E\otimes T+1\otimes E\,, \\
\Delta(U)=U\otimes U\,, & \Delta(F)=F\otimes 1+U\otimes F\,,
\end{array}
\end{equation}
and the duality conditions look like
\begin{equation}
<U\,,T>=R\,,\ \ \ <1\,,T>=<U\,,1>=<F\,,E>={\bf 1}\,. \label{3}
\end{equation}
In~\cite{Vl2} it is shown how an antipode can be introduced into this
bialgebra thus making it a Hopf algebra.

One sees that commutation relations between the $E_i$\,-generators,
as well as between $F^i$, are missing in (\ref{1}). In~\cite{Lu} a
method is proposed how to add such relations to (\ref{1}) without
destroying the bialgebra (the Hopf algebra structure also survives).
The method consists in imposing upon $E$ the relations of the $N$-th
order \begin{equation} E_1\ldots E_N\,\omega_{1\ldots N}=0
\label{4} \end{equation} or, in more explicit form, $$
E_{i_1}\ldots E_{i_N}\omega ^{i_1\ldots i_N}=0\,, $$ for every
solution $\omega $ of the following system of
equations:  \begin{equation} {\left[N\atop n\right]}_B\,\omega
_{1\ldots N}=0 \ \ \ \ (n=1,2,\ldots,N-1)\,, \label{5} \end{equation}
where $B\equiv \check{R}\equiv \sigma R$ and
$ {\left[N\atop n\right]}_B $ are
generalized binomial coefficients~\cite{Lu,Ma-bin}, which are
defined as follows. Let some $x$ and $y$ obey
\begin{equation}
x_1y_2=y_1x_2B_{12}\,. \label{5.5}
\end{equation}
Then
\begin{equation}
(x_1+y_1)\ldots(x_N+y_N)=\sum_{n=0}^{N}y_1\ldots y_nx_{n+1}
\ldots x_N\,{\left[N\atop n\right]}_B\,. \label{6}
\end{equation}

Generalized binomial coefficients enter (\ref{5}) due to the second
relation in the first line of (\ref{1}) which entails
$$ (1\otimes E_1)(E_2\otimes T_2)=(E_1\otimes T_1)(1\otimes E_2)B_{12}
\,. $$ Taking then the coproduct of (\ref{4}) and identifying
$\Delta(E)$ with $y+x$ we can use (\ref{5.5}),(\ref{6}) and
come to the condition (\ref{5}).

It can be easily shown that $\theta =\omega ^t$ (i.e.
$\theta _{i\ldots j}=\omega ^{i\ldots j}$) serves to define the
commutation rules for $F^i$:
\begin{equation}
\theta _{1\ldots N}F_1\ldots F_N=0\,. \label{6.5}
\end{equation}
This completes the definition of our bialgebra (actually, a Hopf
algebra).

The above construction has been illustrated in~\cite{Lu} for the
case of diagonal $R$. Here we choose to consider the simplest
non-diagonal $R$-matrix
\begin{equation}
R=q^\nu \left( \begin{array}{cccc}
q&\cdot &\cdot&\cdot \\ \cdot&1&q-q^{-1}&\cdot \\ \cdot&\cdot &1
&\cdot \\ \cdot&\cdot&\cdot&q  \end{array} \right)\,.  \label{7}
\end{equation}
Dots denote zeros, and a
factor $q^\nu $ is needed to achieve such a normalization
of the $R$\,- or $B$\,-matrix for which the solutions of (\ref{5}) may
exist. Below we give several (not all) explicit solutions for
$\omega $ at $N=2,3,4$ with $\nu $ equal to $1,0,-1/3$, respectively.

Taking into account the explicit form of generalized binomial
coefficients~\cite{Lu,Ma-bin} results in the following equations for
$\omega $ (as usual, $B_i\equiv B_{i\,i+1} $): \\
$N=2\,,\nu =1$:
\begin{equation}
({\bf 1}+B)\,\omega =0\,; \label{8}
\end{equation}
$N=3\,,\nu =0$:
\begin{equation}
\left\{ \begin{array}{l}
({\bf 1}+B_1+B_1B_2)\,\omega =0\,, \\
({\bf 1}+B_2+B_2B_1)\,\omega =0\,;
\end{array} \right. \label{9}
\end{equation}
$N=4\,,\nu =-1/3$:
\begin{equation}
\left\{ \begin{array}{l}
({\bf 1}+B_1+B_1B_2+B_1B_2B_3)\,\omega =0\,, \\
({\bf 1}+B_2+B_2B_1+B_2B_3+B_2B_1B_3+B_2B_1B_3B_2)\,\omega =0\,, \\
({\bf 1}+B_3+B_3B_2+B_3B_2B_1)\,\omega =0\,.
\end{array} \right. \label{10}
\end{equation}
For $N=2$ a solution (tensor $\omega ^{ij}$)
depends on a single parameter $\rho $\,:
$$ \omega ^{11}=\omega ^{22}=0\,,\ \ \omega ^{21}=\rho \,,\ \
\omega ^{12}=-q\rho \,. $$
For $N=3$ it is 2-parametric,
$$ \omega ^{111}=\omega ^{222}=0\,,\ \ \omega ^{211}=\rho \,,\ \
\omega ^{112}=q\rho \,,\ \ \omega ^{221}=\sigma\,, $$
$$ \omega ^{121}=-(1+q)\rho \,,
 \ \ \omega ^{122}=q\sigma \,,\ \
\omega ^{212}=-(1+q)\sigma \,, $$
and for $N=4$ it has 3 parameters:
$$ \omega ^{1111}=\omega ^{2222}=0\,,\ \ \omega ^{2221}=\rho \,,\ \
\omega ^{1222}=-q\rho \,,\ \ \omega ^{2111}=\sigma \,,\ \
\omega ^{1112}=-q\sigma \,, $$
$$ \omega ^{2212}=-q^{-1/3}(1+q^{2/3}+q^{4/3})\rho \,,\ \
\omega ^{2122}=(1+q^{2/3}+q^{4/3})\rho \,,\ \ $$
$$ \omega ^{1211}=-q^{-1/3}(1+q^{2/3}+q^{4/3})\sigma \,,\ \
\omega ^{1121}=(1+q^{2/3}+q^{4/3})\sigma \,, $$
$$ \omega ^{2211}=-\omega ^{1122}=(1+q^{2/3})\varphi \,,\ \
\omega ^{1221}=\omega ^{2112}=(-q^{-2/3}+q^{4/3})\varphi \,, $$
$$ \omega ^{2121}=-(q+2q^{1/3}+q^{5/3})\varphi \,,\ \
\omega ^{1212}=(q^{-1}+q^{-1/3}+2q^{1/3})\varphi \,. $$
These solutions produce the following commutation relations. \\
$N=2$:
\begin{equation}
E_2E_1-qE_1E_2=0\,. \label{11}
\end{equation}
$N=3$:
\begin{equation}
\begin{array}{l} E_2E_1E_1-(1+q)E_1E_2E_1+qE_1E_1E_2=0\,, \\
E_2E_2E_1-(1+q)E_2E_1E_2+qE_1E_2E_2=0\,.
\end{array} \label{12}
\end{equation}
$N=4$:
\begin{equation}
\begin{array}{l}
E_2E_1E_1E_1
+(1+q^{2/3}+q^{4/3})(E_1E_1E_2E_1-q^{-1/3}E_1E_2E_1E_1)
-qE_1E_1E_1E_2=0\,, \\
(q^{1/3}+q^{-1/3})(E_2E_2E_1E_1-E_1E_1E_2E_2)
+(q-q^{-1})(E_1E_2E_2E_1+E_2E_1E_1E_2) \\
\ \ \ \ \ \ \ \ \ \ \ \ \ \ \ \ \ \ \ +(2+q^{-2/3}+q^{-4/3})E_1E_2E_1E_2
-(2+q^{2/3}+q^{4/3})E_2E_1E_2E_1=0\,, \\
E_2E_2E_2E_1+(1+q^{2/3}+q^{4/3})(E_2E_1E_2E_2-q^{-1/3}E_2E_2E_1E_2)
-qE_1E_2E_2E_2=0\,.
\end{array} \label{13}
\end{equation}
The corresponding equations for $F^i$ are the same, with lower indices
replaced by upper ones.

To look more closely at the structure of the resulting algebra
let us arrange its generators as follows:
\begin{equation}
T=\left( \begin{array}{cc} a&0\\b&c \end{array} \right)\,,\ \
U=\left( \begin{array}{cc} x&y\\0&z \end{array} \right)\,,\ \
E=(d\ e)\,,\ \ F =\left( \begin{array}{c} f \\ g \end{array}
\right)
\end{equation}
(zeros are to make the duality relations (\ref{3}) nondegenerate).
The following commutation relations are of
special interest for us: \\
$N=2,3,4$:
\begin{equation} eb-q^\nu
be=q^\nu (q-q^{-1})cd\,,\ \ db-q^{1+\nu }bd=0 \label{15}
\end{equation}
(these stem from (\ref{1})) and, in addition, \\
$N=2$:
\begin{equation}
ed-qde=0\,,  \label{16}
\end{equation}
$N=3$:
\begin{equation}
eed-(1+q)ede+qdee=0\,, \label{17}
\end{equation}
$N=4$:
\begin{equation}
eeed+(1+q^{2/3}+q^{4/3})(edee-q^{-1/3}eede)-qdeee=0\,, \label{18}
\end{equation}
which are the last equations in (\ref{11})-(\ref{13}). Remarkably,
all the other relations in (\ref{12}),(\ref{13}) prove to be dependent
on (\ref{15}), (\ref{17}) and (\ref{18}).

Moreover, it is not difficult to show that our Hopf algebras for
$N=2,3,4$ are isomorphic to $U_qsl(2), U_qsp(2)$ and $U_qg_2$,
respectively. In other words, the $N$-th power in the relations
(\ref{4}) leads to the $q$-deformed universal enveloping algebra $U_qg$
with $g$ defined by the Cartan matrix \begin{equation} \left(
\begin{array}{cc} 2&-1\\-N+1&2 \end{array} \right)\,.  \label{19}
\end{equation} The isomorphism is given by the following relations (we
omit normalization factors which can be easily restored in each
case):  \begin{equation} \begin{array}{llll} \overline{c}b\sim
X_{1}^{+} & a\sim k_{1}^{2}k_{2}^{2} & y\overline{z}\sim X_{1}^{-} &
\overline{x}\sim {k_1'}^2{k_2'}^2 \\ e\sim X_{2}^{+} &
c\sim k_{2}^{2} & g\sim X_{2}^{-} & \overline{z}\sim {k_2'}^2
\end{array} \label{20}
\end{equation}
where $\overline{c}, \overline{x}, \overline{z}$ are inverses
of $c, x, z$, and $\{k_i,X_{i}^{\pm}\}$ the standard Drinfeld-Jimbo
generators of $U_qg$. The second pair $k_1',k_2'$ of Cartan generators
is to be identified with $k_1,k_2$. Second equations in
(\ref{15}), as well as (\ref{16})-(\ref{18}), are associated
(by means of the first equations in (\ref{15})) with
$q$-deformed Serre relations produced by the Cartan matrix (\ref{19}).
An intriguing question whether this correspondence will remain for
$N>4$ is now under investigation.

\vspace{.5cm}

The second topic of the present note concerns quite another procedure
for constructing Hopf algebras, advocated by Manin~\cite{Manin}
in the case of $q$-plane. Recall~\cite{Ma-rev} that the noncommutative
algebras of the coordinate functions on the quantum planes, both
bosonic and fermionic, are respected by coaction of $SL_q(2)$.
Moreover, the requirement of such covariance unambiguously fixes the
bialgebra structure of $SL_q(2)$. Let us see now that, quite
analogously, the bialgebra structure of $sl_q(2)$ is unambiguously
determined if we require the $q$-oscillator algebra to be its covariant
comodule.

The idea of
this possibility has been prompted by Kulish's paper~\cite{Kulish}
where the $sl_q(2)$-comodule structure of the $q$-oscillator has been
discovered.  Let us invert the problem and consider the
$q$-oscillator-type algebra \begin{equation} AB-q^2BA=1 \label{21}
\end{equation}
as a would-be covariant comodule under the coaction of some
bialgebra. Namely, let us require that the commutation relations
(\ref{21}) are respected
by the coaction
\begin{equation} \begin{array}{l} A\rightarrow z\otimes A+x\otimes 1
\,, \\ B\rightarrow z\otimes B+y\otimes 1
\end{array} \label{22}
\end{equation}
of some bialgebra ${\cal A}$ with generators $x,y,z$. Then the
(coassociative) coalgebra structure of ${\cal A}$ is unambiguously
fixed by
\begin{equation} \begin{array}{l} \Delta(z)=z\otimes z\,, \\
\Delta(x)=x\otimes 1+z\otimes x\,, \\
\Delta(y)=y\otimes 1+z\otimes y\,,
\end{array} \label{23}
\end{equation}
and its algebra structure by
\begin{equation} \begin{array}{l} xz=q^2zx\,, \\ zy=q^2yz \,, \\
q^2yx-xy=z^2-1\,. \end{array} \label{24}
\end{equation}
Using the substitution
\begin{equation}
z=q^{-H}\,, \ \ x=\sqrt{q-q^{-1}}q^{-H/2}X^+\,, \ \
y=\sqrt{q-q^{-1}}X^-q^{-H/2} \label{25}
\end{equation}
we come to a conclusion that ${\cal A}$ is precisely $sl_q(2)$\,.
Really, equations (\ref{23}),(\ref{24}) imply
$$ \Delta(X^{\pm})=X^{\pm}\otimes q^{H/2}+q^{-H/2}\otimes X^{\pm}
\,,\ \ \ \Delta(H)=H\otimes 1+1\otimes H\,, $$
$$ [H\,,X^{\pm}]=\pm2X^{\pm}\,,\ \ \ [X^+\,,X^-]=
\frac{q^H-q^{-H}}{q-q^{-1}}\,. $$
This nonstandard approach to deriving $sl_q(2)$ may shed some light on
the role of $q$-oscillator in the general context of quantum groups.

\vspace{.5cm}

{\small This work was supported in part by International Science
Foundation (grant RFF-300) and by Russian Basic Research Foundation
(grant 95-02-05679).

I acnowledge helpful discussions with A.Isaev, P.Kulish, V.Lyakhovsky,
O.Ogievetsky, P.Pyatov, and V.Tolstoy.}

\end{document}